\documentclass[12pt]{amsart}

\usepackage{amsmath}
\usepackage{amsthm}
\usepackage{stmaryrd}
\usepackage{amsfonts,mathrsfs}

\title[$L^2$ Extension with gain]{Analytic inversion of adjunction:\\ $L^2$ Extension Theorems with Gain}
\author{Jeffery D. McNeal$^{\dagger}$}
\thanks{${\dagger}$Partially supported by NSF grant DMS-0300510}
\author{Dror Varolin$^{\ddagger}$}
\thanks{${\ddagger}$Partially supported by NSF grant DMS-0400909}

\newcommand{\noi}{\noindent}

\newcommand{\ms}{\medskip}

\newcommand{\cc}{{\mathcal C}}
\newcommand{\cd}{{\mathcal D}}

\newcommand{\fA}{{\mathfrak A}}

\newcommand{\fD}{{\mathfrak D}}

\newcommand{\vp}{\varphi} 
\newcommand{\ve}{\varepsilon}

\newcommand{\C}{{\mathbb C}}

\newcommand{\R}{{\mathbb R}}

\newcommand{\di}{\partial}
\newcommand{\dbar}{\bar \partial}

\newcommand{\re}{{\rm Re\ }}

\newcommand{\relcomp}{\subset \subset}

\newcommand{\ii}{\sqrt{-1}}

\newcommand{\tensor}{\otimes}

\begin{document}
\maketitle

\theoremstyle{plain}
\newtheorem{thm}{\sc Theorem}
\newtheorem{lem}{\sc Lemma}[section]
\newtheorem{d-thm}{\sc Theorem}[section]
\newtheorem{o-thm}{\sc Theorem}[subsection]
\newtheorem{prop}[lem]{\sc Proposition}
\newtheorem{cor}[lem]{\sc Corollary}

\theoremstyle{definition}
\newtheorem{conj}[lem]{\sc Conjecture}
\newtheorem{defn}[lem]{\sc Definition}
\newtheorem{qn}[lem]{\sc Question}

\theoremstyle{definition}
\newtheorem{ex}{\sc Example}

\theoremstyle{remark}
\newtheorem*{rmk}{\sc Remark}
\newtheorem*{ack}{\sc Acknowledgment}

\section*{Introduction}

\setcounter{section}{1}

The goal of this paper is to establish new results on $L^2$-extension of holomorphic forms 
of top degree with values in a holomorphic line bundle, from an $n$-dimensional submanifold to an $(n+1)$-dimensional manifold.  Motivated by the use of extension results
\begin{enumerate}
\item[(i)] in obtaining estimates for Bergman kernels, and 
\item[(ii)] in the study of multiplier ideals in algebraic geometry, 
\end{enumerate}
we consider the situation where the $L^2$-norms on the submanifold are different from those on the ambient manifold.

We make precise the setting for our work.  Let $X$ be a K\"ahler manifold of complex dimension $n+1$.  Assume there exists a holomorphic function $w$ on $X$, such that $\sup _X |w| = 1$ and, on the set $Z = \{ w=0\}$, $dw$ is never zero.  Thus $Z$ is a smooth hypersurface in $X$.  Assume there exists an analytic subvariety $V \subset X$ such that $Z \not \subset V$ and $X-V$ is Stein.  Thus there are relatively compact subsets $\Omega _j \relcomp X-V$ such that
$$
\Omega _j \relcomp \Omega _{j+1} \quad \text{and} \quad \bigcup _j \Omega _j = X-V.
$$
These hypotheses are satisfied, for example, if $X$ is a Stein manifold, or if $X$ is a holomorphic family of projective algebraic manifolds fibered over the unit disk.  

We assume that $X$ carries a holomorphic line bundle $H \to X$ together with a singular Hermitian metric, given locally by functions $e^{-\kappa}$, whose curvature current is $\ii \di \dbar \kappa.$  Recall that a Hermitian metric is singular if with respect to smooth local trivializations, the functions $\kappa$ are locally integrable.  Thus the curvature of a singular Hermitian metric is a well defined current. We mention that, by this odd but standard terminology (see, for example, \cite{l-04}), smooth Hermitian metrics are singular. Let $R : X \to [-\infty , \infty ]$ be a locally integrable function such that for any local representative $e^{-\kappa}$ of our singular Hermitian mertic over an open set $U$, the function $R+\kappa$ is not identically $+\infty$ or $-\infty$ on $Z\cap U$.  We always assume that 
\begin{equation}\label{i1}
\ii \di \dbar ( \kappa + R + \log |w|^2 ) \ge 0 
\end{equation}
in the sense of currents.

Let $f \in H^0 (Z, H\tensor K_Z)$ be a holomorphic $n$-form on $Z$ with values in $H$.  Following Siu \cite{s-98}, we say that $F \in H^0(X, H\tensor K_X)$ is an extension of $f$ if 
$$
F|Z = f\wedge dw.
$$
We consider the $L^2$-norms
\begin{equation}\label{i2}
N_F= \frac{1}{2\pi} \int _X \frac{|F|^2e^{-\kappa}}{|w|^2 g\left (\log \frac{e} {|w|^2} \right )} \quad \text{and} \quad 
\nu _f= \int _Z |f|^2 e^{-(R+\kappa)}
\end{equation}
for sections of $H\tensor K_X$ and $(H|Z)\tensor K_Z$ respectively.  (Note that if $\sigma$ is a canonical section, then $|\sigma|^2$ is a measure, and can be integrated without reference to a volume form.)

We will formulate conditions on classes of functions $g$ and $R$ that guarantee the existence of an extension $F$ of $f$ satisfying the estimate 
$$
N_F\leq C\nu _f,
$$ 
for any choices of $g$ and $R$ in the conditioned classes, where $C$ is a constant independent of the section $f$. We have written the norms involved in the form \eqref{i2} in order to easily compare our result with previous results. However, 
it is also true that the function $g$, in the norm $N_F$, plays a different role in this problem than that played by $e^{-R}$, in the norm $\nu _f$. Thus we separate these roles by the asymmetric notation in  \eqref{i2}.

\begin{defn}\label{denom} The class $\fD$ consists of nonnegative functions with the following three properties.
\begin{enumerate} 
\item[(i)] Each $g \in \fD$ is continuous and increasing.

\item[(ii)] For each $g \in \fD$ the improper integral
$$
C(g) := \int_1^\infty \frac{dt}{g(t)}
$$
is finite.
\end{enumerate}
For $\delta >0$, set
$$G_\delta(x)=\frac 1{1+\delta}\left( 1+\frac\delta{C(g)}\int_1^x \frac{dt}{g(t)}\right),$$
and note that this function takes values in $(0,1]$. Let
$$
h_{\delta} (x) := \int _1 ^x \frac{1-G_{\delta}(y)}{G_{\delta} (y)} dy.
$$
\begin{enumerate}
\item[(iii)] For each $g \in \fD$ there exists a constant $\delta >0$ such that

$$
K_{\delta}(g) := \sup _{x \ge 1} \frac{x+h_{\delta}(x)}{g(x)}
$$
is finite.
\end{enumerate}
\end{defn}

The conditions on $R$ depend on the choice of the function $g$, and we state these conditions
in the hypotheses of the following Theorem, which is the main result of this paper.

\begin{thm}\label{main}
Let $g$ be a function in $\fD$.  Suppose $R$ is a function such that for all $\gamma > 1$, and $\ve >0$ sufficiently small (depending on $\gamma - 1$),  
\begin{subequations}\label{berg}
\begin{align}
& \alpha - g^{-1} \left ( e^{-R}g(\alpha )\right) \ \text{is subharmonic, and}\\
& g^{-1}\left ( e^{-R}g(1-\log |w|^2) \right ) \ge 1,
\end{align}
\end{subequations} 
where $\alpha = \gamma - \log (|w|^2 + \ve^2).$  Then for every holomorphic $n$-form $f$ with values in $H$ such that 
$$
\int _Z |f|^2 e^{-(R+\kappa)} < +\infty,
$$
there is a holomorphic $(n+1)$-form $F$ with values in $H$ such that 
$$
F|Z = f \wedge dw
$$ 
and 
\begin{equation}\label{mainest}
\frac{1}{2\pi} \int _X \frac{|F|^2 e^{-\kappa}}{|w|^2 g\left ( \log \frac{e}{|w|^2}\right )} \le 4 \left ( K_{\delta}(g) + \frac{1+\delta}{\delta} C(g)\right ) \int _Z |f|^2 e^{-(R+\kappa)}.
\end{equation}
\end{thm}

\begin{rmk}
For any $\lambda > 0$, $K_{\delta} (\lambda g) = \lambda ^{-1} K_{\delta} (g)$ and $C(\lambda g) = \lambda ^{-1} C(g)$.  We thus assume from here on that $\int _1 ^{\infty} \frac{dt}{g(t)} = 1.$  With this particular normalization, we have
\begin{equation*}
\frac{1}{2\pi} \int _{|w|<1} \frac{\ii dw \wedge d\bar w}{|w|^2g\left (\log \frac{e}{|w|^2}\right )}=2 \int _0^1 \frac{dr}{r g\left ( \log \frac{e}{r^2}\right )}= 1\\
\end{equation*}
\end{rmk}

Observe that the function $R \equiv 0$ clearly satisfies the hypotheses \eqref{berg}.  In this situation, we say that Theorem \ref{main} is a theorem on {\it Analytic Inversion of Adjunction}.  For more on Inversion of Adjunction, we refer the reader to the book \cite{l-04} of R. Lazarsfeld or the notes \cite{k-97} of J. Koll\' ar.

The presence of the function $R$ allows us to extend more singular sections on $Z$ to $X$.  Doing so then allows one to get improved estimates for the Bergman kernel.  The type of function $R$ that should be used for good estimates depends on the geometry of the space $X$ and on the way in which $Z$ sits inside $X$.

The seemingly complicated condition (iii) in Definition \ref{denom} is not complicated in practice.  We illustrate this point with examples. Let $g : [1,\infty] \to [0,\infty]$ be one of the functions
\begin{enumerate}
\item [(fn1)]  $$g(x) = s^{-1} e^{s(x-1)}, \quad s \in (0,1]$$

\item [(fn2)]  $$g(x) = x^2$$

\item [(fn3)]  $$g(x) = s^{-1} x^{1+s}, \quad s \in (0,1]$$

\item [(fn4)]  $$g(x) = s^{-1} xL_1(x)L_2(x) \cdots L_{N-2}(x) (L_{N-1}(x))^{1+s}, \quad s\in (0,1],$$
where 
$$
E_j = \exp ^{(j)}(1) \quad \text{and} \quad L_j (x) = \log ^{(j)}
(E_{j}x).
$$
\end{enumerate}
As we shall show in Section \ref{proof-section}, the functions (fn1)-(fn4) all satisfy the conditions in Definition \ref{denom}, as well as the normalization
$$
\int _1 ^{\infty} \frac{dt}{g(t)} =1.
$$ 
Studying these example denominators led us to Theorem \ref{main}.  The cases of Theorem \ref{main} to which they correspond appear ripe for application, and in anticipation we compute bounds on the constants for these examples later in the paper.  The special cases of Theorem \ref{main} corresponding to $g \in {\rm (fn1)}, {\rm (fn2)}$ are already known through previous work.  (This is in part why we have separated (fn2) from (fn3).)  We now recall some of the known results, and compare them with results following from Theorem \ref{main}.

\ms

\noi {\bf The Ohsawa-Takegoshi Extension Theorem:}  If in Theorem \ref{main} we take $R=0$ and $g(x) = e^{x-1}$, we obtain the celebrated Ohsawa-Takegoshi theorem \cite{ot-87}, as phrased by Siu in \cite{s-02}:  For any section $f$ of $K_Z \tensor H |Z $ there exists a section $F$ of $K_X \tensor H$ such that 
$$
F|Z = f \wedge dw \quad \text{and} \quad \int _X |F|^2 e^{-\kappa} \le C \int _Z |f|^2 e^{-\kappa}.
$$
We mention that Berndtsson \cite{b-96} showed that $C$ could be taken to be $8\pi$ (the constant
is written as $4\pi$ in \cite{b-96} because of a slightly different normalization of $w$-area form there). 

\ms

\noi {\bf Ohsawa's negligible weights theorem:}
In the case where $g(x) = e^{x-1}$, we obtain an extension $F$ of $f$ satisfying the estimate
$$
\int _X |F|^2 e^{-\kappa} \le C \int _Z |f|^2 e^{-(R+\kappa)}
$$
under the hypotheses
\begin{subequations}\label{conditions}
\begin{align}
& R + \log |w|^2 \le 0,\\
& \ii \di \dbar R \ge 0,\ \text{and} \\
& \ii \di \dbar (R +\kappa + \log |w|^2) \ge 0.
\end{align}
\end{subequations}
This result was essentially proved by Ohsawa \cite{o-95} under the stronger assumption that, in addition to (\ref{conditions}a), $R$ and $\kappa$ are both plurisubharmonic.  (Ohsawa considered functions on domains in $\C ^n$, which can be identified with canonical sections in that setting.)  

In our situation $\kappa$ need not be plurisubharmonic.  Said another way, we can add $L$ to $R$ and subtract $L$ from $\kappa$ to obtain an estimate
$$
\int _X |F|^2 e^{L-\kappa} \le C \int _Z |f|^2 e^{-(R+\kappa)}
$$
for an extension $F$ of $f$, under the hypotheses
\begin{eqnarray*}
&& L+R + \log |w|^2 \le 0,\\
&& \ii \di \dbar (L+R) \ge 0,\ \text{and} \\
&& \ii \di \dbar (R +\kappa + \log |w|^2) \ge 0.\\
\end{eqnarray*}

\noi {\bf Demailly's logarithmic extension theorem:}
In the case $R=0$ and $g= x^2$, we obtain the following result of Demailly \cite[Theorem 12.6]{d-00}.

\begin{d-thm}[Demailly]\label{demailly-log}
For any holomorphic section $f$ of $K_Z\tensor H|Z$ such that
$$
\int _Z |f|^2 e^{-\kappa} < +\infty
$$
there exists a section $F$ of $K_X \tensor H$ such that 
$$
F|Z=f\wedge dw \quad \text{and} \quad \int _X \frac{|F|^2 e^{-\kappa}}{|w|^2 \left ( \log \frac{e}{|w|^2} \right )^2}
\le C \int _Z |f|^2 e^{-\kappa}.
$$
\end{d-thm}

Demailly proved Theorem \ref{demailly-log} by somewhat different methods.  Demailly's proof is more complicated than our proof of Theorem \ref{main} below, perhaps in part because he considers the more general case in which $w$ is a section of some vector bundle $E \to X$, and the zero locus $Z$ is only generically smooth.

Demailly points out that the trivial estimate 
$$
s \log \frac{e}{x} \le \frac{1}{x^s}
$$ 
implies that the extension $F$ satisfies 
$$
\int _X \frac{|F|^2 e^{-\kappa}}{|w|^{2-2s}} \le \frac{C}{s^2} \int _Z |f|^2 e^{-\kappa}.
$$
As we will show in section 3, Theorem \ref{main}, applied to the function $g(x) = s^{-1} e^{s(x-1)}$, gives an extension $\tilde F$ with the better estimate
\begin{equation}\label{alg-est}
\int _X \frac{|\tilde F|^2 e^{-\kappa}}{|w|^{2-2s}} \le \frac{C }{s} \int _Z |f|^2 e^{-\kappa}.
\end{equation}
Moreover, for any extension $\hat F$ of $f$ such that 
$$
\int _X \frac{|\hat F|^2 e^{-\kappa}}{|w|^{2-2s}} < +\infty
$$
the function 
$$
s \mapsto \int _X \frac{|\hat F|^2 e^{-\kappa}}{|w|^{2-2s}}
$$
extends to a meromorphic function with a simple pole at $0$.
%The residue of this function may be computed as follows.  By approximation we may assume that 
%\begin{enumerate}
%\item[(i)] $w$ is smooth up to the boundary, 
%\item[(ii)] for some $\ve > 0$, $|dw|$ is nowehere $0$ on the set $Z_{\ve} := \{ |w| < \ve \}$, and 
%\item[(iii)] On the set $Z_{\ve}$, $|\hat F|^2 e^{-\kappa}$ is smooth in the variable $w$.
%\end{enumerate}
%Our assumptions imply that, perhaps after shrinking, $Z_{\ve} = Z \times \{ |w| < \ve\}$.  We can then write 
%$$
%|\hat F|^2 e^{-\kappa} (z,w) = \left (|f|^2e^{-\kappa (z,0)} + w H_1 (z) + \bar w H_2 (z) + O(|w|^2) \right )  \ii dw \wedge d\bar w
%$$
%for $z \in Z$ and $|w| < \ve$.  Thus 
%\begin{eqnarray*}
%s \int _X \frac{|\hat F|^2 e^{-\kappa}}{|w|^{2-2s}} &=& s \int _{Z_{\ve}}  \frac{|\hat F|^2 e^{-\kappa}}{|w|^{2-2s}} + O(s) = 2\pi \ve ^{2s}\int _Z |f|^2 e^{-\kappa} + O(s) 
%\end{eqnarray*}
%Taking the limit as $s \to 0$ shows that the residue is 
%$$
%2\pi \int _Z |f|^2 e^{-\kappa}.
%$$
Thus the order of growth of the coefficient in the bound obtained for $\tilde F$ is the best possible.

\ms

\noi {\bf Logarithms in the estimates:}
We point out that the cases when $g \in {\rm (fn3)}, {\rm (fn4)}$ in Theorem \ref{main} give significantly
refined estimates on the extension $F$ as compared with \eqref{alg-est}. For example, if 
$g \in {\rm (fn3)}$ and $R=0$, Theorem \ref{main} gives an extension satisfying the estimate
$$
\int _X \frac{|F|^2e^{-\kappa}}{|w|^2 \left ( \log \left ( \frac{e}{|w|^2}\right ) \right ) ^{1+s}} \le \left ( \frac{1}{s} +1 \right ) 8 \pi  \int _Z |f|^2 e^{-\kappa} \quad \text{as} \ s \searrow 0.
$$
The ability to hold the power of $|w|$ at 2 in the denominator on the left-hand side in this estimate, while pushing the integrand on the left-hand side towards non-integrability by letting $s\rightarrow 0$, opens new possibilities of obtaining extensions which satisfy additional side properties. 

\ms

We give the proof of Theorem \ref{main} in the next section. In section 3, we consider the special
cases $g\in\text{(fn1)-(fn4)}$.

\section{The general framework}\label{gen}

\subsection{Twisted Bochner-Kodaira Identity and Basic Estimate}

We begin by reviewing the now standard twisted Bochner-Kodaira identity.

Let $\Omega$ be a smoothly bounded pseudoconvex domain in our K\" ahler manifold $X$.  Consider the usual Bochner-Kodaira identity for $H$-valued $(n+1,1)$-forms $u$. If $u$ is such a form that is also in Dom $(\dbar^*_{\vp})$ and has components in $C^\infty(\overline\Omega)$,  in which case we shall write $u\in\cd$, then:
\begin{eqnarray}\label{bk}
&& \int _{\Omega} \left | \dbar ^* _{\vp} u\right |^2 e^{-\vp} + \int _{\Omega} \left | \dbar u\right |^2 e^{-\vp} = \int _{\Omega} \ii \di \dbar \vp (u,u) e^{-\vp}\\
\nonumber && \qquad + \int _{\Omega} \left | \overline \nabla u\right |^2 e^{-\vp} +
\int _{\di \Omega}  \ii \di \dbar \rho (u,u) e^{-\vp}
\end{eqnarray}
Note that the last two terms of the right-hand side
are non-negative, the last since $\Omega$ is pseudoconvex. We mention that $u\in\cd$ forces a boundary condition on $u$ which allows \eqref{bk} to take the stated form.

Next we pass the the twisted estimates.  The main philosophical idea is to consider a twist of the original metric for $H$.  That is to say, we consider a metric $e^{-\psi}$ for $H$.  For any such metric, there exists a positive function $\tau$ such that
$$
e^{-\vp} = \tau e^{-\psi}.
$$
Now, 
$$
\di \dbar \vp = \di \dbar \psi - \tau ^{-1} \di \dbar \tau + \tau ^{-2} \di \tau \wedge \dbar \tau.
$$
Also, using the formula
$$
\dbar ^* _{\vp} u = - \sum _j e^{\vp} \frac{\di}{\di z^j} (e^{-\vp}u_{\bar j}),
$$
where locally $u = \sum u_{\bar j} d\bar z ^j$ with $u_{\bar j}$ canonical sections, we have the formula
$$
\dbar ^*_{\vp} u =-\tau ^{-1} \di \tau (u) + \dbar ^* _{\psi} u.
$$
Substitution of these identities into (\ref{bk}) then gives, after a certain amount of educated manipulation, the so-called twisted Bochner-Kodaira identity: if $u\in\cd$ then
\begin{eqnarray}\label{tbk}
&& \int _{\Omega} \tau \left | \dbar ^* _{\psi} u\right |^2
e^{-\psi} + \int _{\Omega} \left | \dbar u\right |^2 e^{-\vp} \\
\nonumber && \qquad = \int _{\Omega}\left (\tau \ii \di \dbar \psi -  \ii \di \dbar \tau \right ) (u,u) e^{-\psi}+ 2 \re \int _{\Omega} \di \tau (u) \overline{\dbar ^* _{\psi} u} e^{-\psi}\\
\nonumber && \qquad \qquad + \int _{\Omega} \tau \left | \overline \nabla u\right |^2 e^{-\psi} + \int _{\di \Omega}  \ii \di \dbar \rho (u,u) e^{-\vp}
\end{eqnarray}
This identity was obtained in \cite{ot-87} for forms $u$ with compact support in $\Omega$. For
$u\in\cd$, the identity was obtained independently by \cite{b-96}, \cite{mc-96}, and \cite{s-96}.

By using the positivity of the terms on the last line of (\ref{tbk}) and applying the Cauchy-Schwarz inequality to the third term in the second line, we obtain this lemma.

\begin{lem}\label{mcneal-lem}
For any $(n+1,1)$-form $u$ in the domain of the adjoint of $\dbar_\psi ^*$, the following inequality holds.
\begin{eqnarray}\label{mcneal-id}
&&\int _{\Omega} (\tau + A) \left | \dbar ^* _{\psi} u\right |^2 e^{-\psi}  + \int _{\Omega} \tau \left | \dbar u \right | ^2
e^{-\psi} \\
\nonumber && \qquad \ge \int _{\Omega} \left ( \tau \ii \di \dbar \psi -
\ii \di \dbar \tau  - \frac{1}{A} \ii \di \tau \wedge \dbar \tau \right ) (u,u)  e^{-\psi}.
\end{eqnarray}
\end{lem}
To obtain \ref{mcneal-lem} from \eqref{tbk}, one also needs that $\cd$ is dense in Dom 
$(\dbar_\psi ^*)$ in the graph norm given by the left-hand side of \eqref{tbk}. We refer to \cite{d-00}
for this fact.

\subsection{Choices for $\tau$, $A$ and $\psi$, and an a priori estimate}

Fix a constant $\gamma > 1$.  (Eventually we will let $\gamma \to 1$.)  We define the function $a$ to be
\begin{equation}\label{s=1}
a := g^{-1} \left ( e^{-R}g\left (\gamma - \log (|w|^2+ \ve ^2) \right ) \right ),
\end{equation}
where $g\in\fD$.  We note that, by Condition (\ref{berg}b), $a > 1$ if $\ve > 0$ is sufficiently small.  In view of property (\ref{berg}a),  a straightforward calculation shows that
$$
- \ii \di \dbar a \ge \frac{\ve ^2}{(|w|^2 + \ve ^2)^2} \ii \di w \wedge d\bar w.
$$

Suppose $h : [1,\infty) \to \R $ is a function such that for all $x>1$ and some constant $M > 0$, 
\begin{subequations}\label{h-conditions}
\begin{align}
& x+h(x) \ge 1,\\
& 1+ h'(x) \ge 1,\ \text{ and} \\
& h''(x) < 0.
\end{align}
\end{subequations}
We then take 
\begin{equation}\label{tau-a-choices-s=1}
\tau = (a + h(a))  \quad \text{and} \quad A = \frac{(1+h'(a))^2}{-h''(a)}.
\end{equation}
With these choices, one calculates that 
\begin{equation}\label{tau-A}
-\ii \di \dbar \tau - \frac{\ii}{A} \di \tau \wedge \dbar \tau = (1+h'(a)) (-\ii \di \dbar a) \ge \frac{\ve^2}{(|w|^2+\ve ^2)^2} \ii dw \wedge d\bar w.
\end{equation}

Before proceeding, we show how to choose $h$, and, hence $\tau$ and $A$, from the given function $g\in\fD$.

\begin{lem}\label{ode-lem} If $g\in\fD$ and $\int _1 ^{\infty} \frac{dt}{g(t)} = 1$, then there exists an $h$ satisfying conditions {\rm (\ref{h-conditions}a)-(\ref{h-conditions}c)} and the ODE
\begin{equation}\label{ode}
h''(x) +\frac {\delta}{(1+\delta) g(x)}\left(1+h'(x)\right)^2 =0,\qquad x\geq 1,
\end{equation}
where $\delta$ is a positive number whose existence is guaranteed by Definition \ref{denom}. 
\end{lem}

\begin{proof} For the given $g\in\fD$, let $h=h_\delta$ be the function given by Definition \ref{denom}, i.e.
$$
h(x)= \int _1 ^x \frac{1-G_{\delta}(y)}{G_{\delta} (y)} dy,
$$
where $G_\delta$ is given by
$$
G_\delta(x)=\frac 1{1+\delta}\left( 1+\delta \int_1^x \frac{dt}{g(t)}\right).
$$
It is straightforward to check that $h$ satisfies \eqref{ode}, and that condition (\ref{h-conditions}c) is satisfied. For the remaining conditions, note that $1+h'(x)= \frac 1{G_\delta(x)}$, so (\ref{h-conditions}b) follows from the inequality $G_\delta(x)\le 1$.  Next,
$$
x+h(x) =\int _1 ^x 1+\frac{1-G_{\delta}(y)}{G_{\delta} (y)} dy +1 \ge 1,
$$
since $G_\delta >0$, which shows that (\ref{h-conditions}a) holds.
\end{proof}

Finally, we take 
$$
\psi = \kappa + R + \log |w|^2.
$$
With these choices, substitution into Lemma \ref{mcneal-lem} gives us the following {\it a priori} estimate:
\begin{equation}\label{a-priori-est}
\int _{\Omega} \frac{\ve ^2}{(|w|^2 + \ve^2)^2} | \left < u , d \bar w \right > |^2 e^{-\psi} \le ||T^* u||_\psi^2 + ||Su||_\psi^2 ,
\end{equation}
where 
$$
T\beta  = \dbar \left ( \sqrt{\tau + A} \beta \right ) \quad \text{and} \quad Su = \sqrt{\tau }\left ( \dbar u\right ).
$$

\subsection{A smooth extension and its holomorphic correction}

Since $\Omega$ is Stein, we can extend $f$ to an $H$-valued holomorphic $n$-form $\tilde f$ on $\Omega$.  By extending to a Stein neighborhood of $\Omega$ (which exists by hypothesis) we may also assume that 
$$
\int _{\Omega}\left | \tilde f \wedge dw \right |^2 < +\infty.
$$
Of course, we have no better estimate on this $\tilde f$.  In particular, the estimate would degenerate as $\Omega$ grows.

In order to tame the growth of this extension $\tilde f$, we first modify it to a smooth extension.  To this end, let $\delta >0$ and let $\chi \in \cc ^{\infty} _0 ([0,1))$ be a
cutoff function with values in $[0,1]$ such that
$$
\chi \equiv 1 \ {\rm on} \ [0,\delta] \quad {\rm and} \quad |\chi
'|\le 1+\delta.
$$
We write
$$
\chi _{\ve} := \chi \left ( \frac{|w|^2}{\ve^2} \right ).
$$
We distinguish the $(n+1,1)$-form
$$
\alpha _{\ve} := \dbar \chi _{\ve} \tilde f \wedge dw.
$$
Then one has the estimate
\begin{eqnarray*}
\left|(u,\alpha _{\ve})_\psi\right|^2 &\le & \left ( \int _{\Omega}|\left <
u,\alpha _{\ve}\right >| e^{-\psi} \right )^2 \\
&=& \left ( \int _{\Omega} \left | \left < u, \frac{w}{\ve^2} \chi
' \left ( \frac{|w|^2}{\ve ^2}\right ) \tilde
f \wedge dw \wedge d\bar w \right >\right | e^{-\psi} \right )^2 \\
&\le & \int _{\Omega} \left | \frac{\tilde f \wedge dw}{\ve ^2}
\chi ' \left ( \frac{|w|^2}{\ve ^2}\right )\right |^2
\frac{(|w|^2+\ve ^2)^2}{\ve^{2}} e^{-(\kappa+R)} \\
&& \times \int _{\Omega} |\left < u ,d\bar w\right > |^2
\frac{\ve^{2}}{(|w|^2+\ve ^2)^{2}} e^{-\psi} \\
&\le & \frac{C_{\ve}}{M} \left ( ||T^*u||_\psi^2 + ||Su||_\psi^2 \right )
\end{eqnarray*}
where
$$
C_{\ve} := \frac{4(1+\delta)^2}{\ve ^2} \int _{|w| \le \ve} \left
| \tilde f \wedge dw \right |^2 e^{-(\kappa+R)} .
$$
Note that the factor $|w|^2$ cancels the term $\log|w|^2$ in $\psi$ in the third line above.

By the usual $L^2$-method we obtain the following result.
\begin{d-thm}\label{t-est-thm}
There exists a smooth $(n+1)$-form $\beta _{\ve}$ such that
$$
T \beta _{\ve} = \alpha _{\ve} \quad {\rm and} \quad \int _{\Omega}
|\beta _{\ve}|^2 e^{-\psi} \le \frac{C_{\ve}}{M}.
$$
In particular,
$$
\beta _{\ve} |Z \equiv 0.
$$
\end{d-thm}
\begin{proof} Let $\fA:=\left\{T^*u:u\in\fD\right\}$.
On the set $\fA\cap\text{Null }(S)$, the anti-linear functional
$$
\ell : T^*u\longrightarrow (\alpha_\epsilon, u)_\psi
$$ 
is bounded, by the inequality above. On $\left[\text{Null }(S)\right]^\perp$, the functional is trivially bounded. By extending $\ell$ trivially in the direction orthogonal to the image of $T^*$, we can assume, without increasing the norm, that $\ell$ is defined on the whole Hilbert space.  The Riesz representation theorem then gives a solution to $T \beta _\epsilon = \alpha_\epsilon$ with the stated norm inequality on $\beta_\epsilon$.

It remains only that $\beta _{\ve} |Z \equiv 0$.  But one notices that $\psi$ is at least as singular as $\log |w|^2$, and thus $e^{-\psi}$ is not locally integrable at any point of $Z$.  The desired vanishing of $\beta _{\ve}$ follows.
\end{proof}

The candidate for the extension of $f \wedge dw$ will be the holomorphic section
$$
F = \limsup _{\ve \to 0} \left ( \chi _{\ve} \tilde f \wedge dw - \sqrt{\tau + A} \beta _{\ve}\right ).
$$
As a result, we obtain the following estimates for $F$:
$$
\frac{1}{2\pi} \int _{\Omega} \frac{|F|^2e^{-\kappa}}{|w|^2g(1-\log |w|^2)}\le \limsup _{\ve \to 0} \sup _X \left ( \frac{(\tau +A)e^{R}} {g(1-\log |w|^2)} \right ) \times 4\int _Z |f|^2 e^{-(R+\kappa)}.
$$
Indeed, the first of the two terms in the definition of $F$ is supported on a set whose measure converges to 0 with $\ve$, and thus one has only to estimate the second term using Theorem \ref{t-est-thm}.

Thus it remains only to show that the quantity 
$$
\limsup _{\ve \to 0, \gamma \to 1} \sup _X \left ( \frac{(\tau +A)e^{R}}{ g(1-\log |w|^2)} \right ) = \sup _{a \ge \gamma _o} \frac{\tau +A}{g(a)}.
$$
is finite. However, $\tau$ was chosen according to \eqref{tau-a-choices-s=1}, i.e. $\tau(x)=x+h(x)$, so by Definition \ref{denom} (iii)
$$
\sup_{x\ge 1} \frac{\tau(x)}{g(x)}=K_\delta(g)< +\infty,
$$
and $A$ was also chosen in \eqref{tau-a-choices-s=1}, i.e. $A=\frac{(1+h')^2}{-h''}$, so by Lemma \ref{ode-lem}
$$
\frac{A(x)}{g(x)} = \frac{1+ \delta}{\delta} < +\infty.
$$
This completes the proof of Theorem \ref{main}.

\section{Specific denominators and constants} \label{proof-section}

We now examine the cases of $g\in\text{(fn1)-(fn4)}$.  The purpose of this section is to show how easily the criteria of Definition \ref{denom} can be verified for these examples.  Thus we do not worry too much about obtaining the best estimates.  However, sharper estimates can be obtained by a more careful analysis than we carry out here.

\subsection{Algebraic gain: the case $g(x) = s^{-1}e^{s(x-1)}$}
First one has
$$
\int _1 ^x \frac{dt}{g(t)} = 1 - e^{-s(x-1)}.
$$
Lead by Definition \ref{denom}, we compute that 
$$
G_{\delta} (x) = \frac{(1 + \delta)e^{s(x-1)} - \delta )}{(1+\delta) e^{s(x-1)}}
= 1 - \frac{\delta}{(1+\delta) e^{s(x-1)}}
$$ 
and thus, with $\ve = \frac {\delta}{1+\delta}$, that 
\begin{eqnarray*}
h_{\delta} (x) &=& \int _1 ^x \frac{1-G_{\delta}(y)}{G_{\delta} (y)}dy   = \int _1 ^x \frac{\ve dy}{e^{s(y-1)} - \ve}\\
&\le & (x -1) \frac{\ve}{1-\ve} = (x-1)\delta.
\end{eqnarray*}
It follows that 
$$
\frac{s(x+h_{\delta}(x))}{e^{s(x-1)}} \le \frac{s ((1+\delta) (x-1) + 1)}{e^{s(x-1)}} = \frac{(1+\delta) r + s}{e^r},
$$
where $r = s(x-1)$.  In general, the function $e^{-r}(ar+b)$, $r \ge 0$, achieves its maximum at the point $r = (a-b)/a$.  It follows that 
$$
K_{\delta}(g) \le \frac{1+\delta}{\exp\left ( \frac{1+\delta - s}{1+\delta}\right )} \le 1+\delta.
$$
Applying Theorem \ref{main}, we obtain an extension $F$ of $f$ with the estimate 
$$
\frac{1}{2\pi} \int _X \frac{|F|^2 e^{-\kappa}}{|w|^{2-2s}} \le \frac{(1+\delta)^2}{\delta} \frac{4}{s} \int _Z |f|^2 e^{-\kappa}.
$$
Incidentally, the right hand side is minimized by taking $\delta =1$.

\subsection{Logarithmic denominator:  the case $g(x)=x^2$}

First we have
$$
\int _1 ^x \frac{dt}{g(t)} = 1-\frac{1}{x}.
$$
Then 
$$
G_{\delta} (x) = \frac{(1+ \delta) x - \delta }{(1+\delta)x},
$$
and thus
$$
h_{\delta} (x) = \int _1 ^x \frac{1-G_{\delta}(y)}{G_{\delta}(y)} dy = \int _1 ^x \frac{\delta  dy}{(1+\delta) y-\delta} \le \delta (x-1).
$$
Then 
$$
\frac{x+h_{\delta}(x)}{x^2} = \frac{(1+\delta)}{x} - \frac{\delta}{x^2} \le \frac{(1+\delta )^2 }{4\delta}.
$$
Applying Theorem \ref{main}, we obtain an extension $F$ of $f$ with the estimate
$$
\frac{1}{2\pi} \int _X \frac{|F|^2e^{-\kappa}}{|w|^2 \left ( \log \frac{e}{|w|^2}\right )^2} \le (2+\delta) \frac{1+\delta}{\delta} \int _Z |f|^2 e^{-\kappa}.
$$
The value $\delta = \sqrt 2$ minimizes the right hand side, in which case
$$
 (2+\delta) \frac{1+\delta}{\delta} = 3+2 \sqrt 2.
 $$

\subsection{Logarithmic inversion of adjunction:  the case $g(x) = s^{-1} x^{1+s}$}

This time 
$$
\int _1 ^x \frac{dt}{g(t)} = 1-\frac{1}{x^s}.
$$
Then 
$$
G_{\delta} (x) = \frac{(1+\delta)x^s - \delta }{(1+\delta)x^s} =  1 - \frac{\delta}{(1+\delta)x^s},
$$
and thus with $\ve = \frac{\delta }{1+\delta}$ we have 
$$
h_{\delta}(x) =\int _1 ^x \frac{1-G_{\delta}(y)}{G_{\delta}(y)} = \int _1 ^x \frac{\ve dy}{y^s-\ve} \le (x-1) \frac{\ve}{1-\ve} = \delta (x-1).
$$
It follows that 
$$
\frac{s(x+h_{\delta}(x))}{x^{1+s}} = \frac{s(1+\delta)x - s\delta}{x^{1+s}} \le s(1+\delta).
$$
We obtain from Theorem \ref{main} an extension $F$ of $f$ satisfying the estimate 
$$
\frac{s}{2\pi} \int _X \frac{|F|^2 e^{-\kappa}}{|w|^2 \left ( \log \frac{e}{|w|^2} \right ) ^{1+s}} \le \left ( \frac{(1+\delta s)(1+\delta)}{\delta} \right ) 4 \int _Z |f|^2 e^{-\kappa}.
$$
At this point we can let $\delta \to \infty$.  For example, we can take $\delta = s^{-1/2}$.  We then obtain an extension $F$ of $f$ satisfying the estimate
$$
\frac{s}{2\pi} \int _X \frac{|F|^2 e^{-\kappa}}{|w|^2 \left ( \log \frac{e}{|w|^2} \right ) ^{1+s}} \le 4 ( 1+2\sqrt{s} + s)  \int _Z |f|^2 e^{-\kappa}.
$$

\subsection{Iterated logs: the case $g \in {\rm (fn4)}$}
Recall that
$$
E_j = \exp ^{(j)}(1) \quad \text{and}\quad L_j(x) =\log ^{(j)}(E_jx),
$$
and that 
$$
g(x) = s^{-1} x \left ( \prod _{j=1}^{N-2} L_j (x) \right ) (L_{N-1}(x))^{1+s}.
$$
Then 
$$
\int _1 ^x \frac{dt}{g(t)} = 1 - \frac{1}{(L_{N-1}(x))^s}
$$
and thus 
$$
G_{\delta} (x) = \frac{(1+\delta)(L_{N-1}(x))^s - \delta }{(1+\delta)(L_{N-1}(x))^s} =  1 - \frac{\delta}{(1+\delta)(L_{N-1}(x))^s}.
$$
Putting $\ve = \frac{\delta }{1+\delta}$, we have 
\begin{eqnarray*}
&& h_{\delta}(x) =\int _1 ^x \frac{1-G_{\delta}(y)}{G_{\delta}(y)} = \int _1 ^x \frac{\ve dy}{(L_{N-1}(y))^s-\ve} \\ &&\le (x-1) \frac{\ve}{1-\ve} = \delta (x-1).
\end{eqnarray*}
It follows that 
\begin{equation*}\label{log-bound}
\frac{x+h_{\delta}(x)}{g(x)} = \frac{ s (1+\delta) - \frac{\delta}{x} }{\left (  \left ( \prod _{j=1}^{N-2} L_j (x) \right ) (L_{N-1}(x))^{1+s} \right )^{-1}} \le s(1+\delta).
\end{equation*}
As in the previous paragraph, we obtain an extension $F$ of $f$ satisfying the estimate 
$$
\frac{1}{2\pi} \int _X \frac{|F|^2 e^{-\kappa}}{|w|^2  g \left ( \log \frac{e}{|w|^2} \right ) } \le 4 ( 1+2\sqrt s + s)  \int _Z |f|^2 e^{-\kappa}.
$$

\ms


\begin{thebibliography}{99}

\bibitem[B-96]{b-96} Berndtsson, B., {\it The extension theorem of Ohsawa-Takegoshi and the theorem of Donnelly-Fefferman.}  Ann. Inst. Fourier (Grenoble) 46 (1996), no. 4, 1083--1094. 

\bibitem[D-00]{d-00} Demailly, J.-P., {\it Multiplier ideal sheaves and analytic methods in algebraic geometry.}  School on Vanishing Theorems and Effective Results in Algebraic Geometry (Trieste, 2000), 1--148, ICTP Lect. Notes, 6, Abdus Salam Int. Cent. Theoret. Phys., Trieste, 2001.

\bibitem[K-97]{k-97} Koll\'ar, J., {Singularities of pairs.}  Algebraic geometry---Santa Cruz 1995, 221--287, Proc. Sympos. Pure Math., 62, Part 1, Amer. Math. Soc., Providence, RI, 1997.

\bibitem[L-04]{l-04} Lazarsfeld, R., {\it Positivity in algebraic geometry, I, II}, Springer (2004)

\bibitem[M-93]{m-93} Manivel, L., {\it Un th\' eor\` eme de prolongement $L\sp 2$ de sections holomorphes d'un fibr\' e hermitien.} Math. Z. 212 (1993), no. 1, 107--122.

\bibitem[Mc-96]{mc-96} McNeal, J.D.,
{\it On large values of $L^2$ holomorphic functions}
Math. Res. Let. 3 (1996), 247-259.

\bibitem[O-95]{o-95} Ohsawa, T., {\it On the extension of $L\sp 2$ holomorphic functions. III. Negligible weights.}  Math. Z. 219 (1995), no. 2, 215--225.

\bibitem[OT-87]{ot-87} Ohsawa, T., Takegoshi, K., {\it On the extension of $L^2$ holomorphic functions.}  Math. Z. 195 (1987), no. 2, 197--204. 

\bibitem[S-96]{s-96} Siu, Y.-T., {\it The Fujita conjecture and the extension theorem
of Ohsawa-Takegoshi} Geometric Complex Analysis, Hayama. World Scientific
(1996), 577-592.

\bibitem[S-98]{s-98} Siu, Y.-T., {\it Invariance of plurigenera.}  Invent. Math. 134 (1998), no. 3, 661--673.

\bibitem[S-02]{s-02}Siu, Y.-T., {\it Extension of twisted pluricanonical sections with plurisubharmonic weight and invariance of semipositively twisted plurigenera for manifolds not necessarily of general type.}  Complex geometry. Collection of papers dedicated to Hans Grauert.  Springer-Verlag, Berlin, 2002. (223--277)

\end{thebibliography}
\end{document}